\documentclass{article}

\usepackage{amsthm}

\newtheorem{theorem}{Theorem}[section]

\newtheorem{lemma}[theorem]{Lemma}

\newtheorem{corollary}[theorem]{Corollary}
\newtheorem{conjecture}[theorem]{Conjecture}
\newtheorem{definition}[theorem]{Definition}

\newtheorem{remark}[theorem]{Remark}

\usepackage{xcolor}
\usepackage{comment}
\usepackage{tkz-euclide}
\usepackage{ytableau}
\usepackage{mathtools}

\begin{document}

\author{Galyna Dobrovolska \thanks{E-mail address: galdobr@gmail.com}}
\title{Some remarks on combinatorial wall-crossing}

\maketitle

\begin{abstract}

\noindent We establish a new simple explicit description of the Young diagram at each step of the combinatorial wall-crossing algorithm for the rational Cherednik algebra applied to the trivial representation. In this way we provide a short and explicit proof of a key theorem of P. Dimakis and G. Yue. We also present two conjectures on combinatorial wall-crossing which were found using computer experiments.

\end{abstract}

\section{Introduction}

We study an algorithm for combinatorial wall-crossing for the rational Cherednik algebra. More precisely, let $0<\frac{a^{\prime}_1}{b^{\prime}_1}<...<\frac{a^{\prime}_k}{b^{\prime}_k}<a/b$ be all the rational numbers with denominator at most $n$ between $0$ and $a/b$, which we call "walls," i.e. the terms of the $n$-th Farey sequence which are smaller than $a/b$. Let $M_d$ be the generalized Mullineux involution (see Section 2) on the set of $d$-regular partitions. We will study the properties of the permutation $\widetilde{M}_{a/b}= M_{b^{\prime}_k}^t \circ \dots \circ M_{b^{\prime}_1}^t$ on partitions of $n$, where "t" stands for the transposed partition. We call the involution $M_{b^{\prime}_i}^t$ wall-crossing across the wall $\frac{a^{\prime}_i}{b^{\prime}_i}$).

 The study of wall-crossing functors was initiated by Beilinson and Ginzburg (\cite{BG}). Wall-crossing functors were recently used as a crucial ingredient in the work of Bezrukavnikov and Loseu (\cite{BL}). Furthermore Loseu proved in \cite{L} that combinatorial wall-crossing is given by the Mullineux involution in large positive characteristic. Wall-crossing was also recently studied by Gorsky and Negut (\cite{GN}) and Su, Zhao, and Zhong (\cite{SZZ}). 
 
 The motivation for considering the above permutation $\widetilde{M}_{a/b}$ comes from Bezrukavnikov's combinatorial conjecture (see Section 2) which stems from relations between wall-crossing functors and monodromy of the quantum connection for the Hilbert scheme of points in the plane (cf. \cite{BO}). More precisely, Bezrukavnikov conjectured that certain invariants (corresponding to dimensions of supports of simple representations of the rational Cherednik algebra) of the composition of the above involutions are equal to similar invariants for a composition of involutions which are much simpler, in particular do not involve the Mullineux involution in their definition (see Section 2 for details). 

We note that a similar study to the present one has been undertaken in \cite{DY}. The goal of the present paper is to streamline, clarify, and add new information to the results obtained in \cite{DY}. Our method of proof is different from the method of proof in \cite{DY}, although we make use of similar general constructions. Our methods have the advantage that they provide more explicit information about the intermediate steps of the algorithm and they lead to the sought conclusion faster.

Now we will state the main theorem of this note. 

We will talk about coordinates of any box in any Young diagram in the following way. We will place the Young diagram in the fourth quadrant (i.e. the quadrant with $x>0$ and $y<0$) so that the horizontal and vertical boundaries of the Young diagram are aligned with the corresponding coordinate axes 
The coordinates of a box in the Young diagram will be the Cartesian coordinates of its top left corner in this placement.

\begin{definition} We define an order on pairs of integers by $(x_1,y_1) >_{a,b} (x_2,y_2)$ if $ax_1-(b-a)y_1 > a x_2-(b-a) y_2$ or $ax_1-(b-a)y_1 = a x_2-(b-a) y_2$ and $y_1>y_2$.
\end{definition}

\begin{definition} Let $D^n_{a,b}$ be the union of $n$ smallest boxes with respect to the order $>_{a,b}$. 
\end{definition}

In this note we will prove the following

\begin{theorem} \label{main}
The result of the combinatorial wall-crossing transformation applied to the trivial partition to the left of the wall $a/b$ is $D^n_{a,b}$ i.e. $\widetilde{M}_{a/b}({\rm triv}_n)=D^n_{a,b}$.
\end{theorem} 

The note is organized as follows. In Section 2 we state the preliminaries, in particular describe Bezrukavnikov's conjecture. In Section 3 we find a new description of partitions obtained by applying the above algorithm to the trivial representation, and using this description we prove a statement equivalent to the main theorem of \cite{DY}. In Section 4 we conjecture the effect of the algorithm on the sign representation for a prime value of the size of the partition. Finally, in Section 5 we conjecture the inductive structure of this algorithm by means of adding one good box.

\subsection*{Acknowledgements}
The author is very grateful to Roman Bezrukavnikov for suggesting this problem and also to him and Panagiotis Dimakis and Guangyi Yue for many useful discussions and ideas. The author was supported by an NSF Postdoctoral Fellowship. 

\section{Preliminaries}

For positive integers $a<b$ we define the operation of generalized column regularization ${\rm colreg}_{a,b}$ (which first appeared in \cite{DY}) on partitions. The operation consists of sliding all boxes on any ladder which is the collection of points with integer coordinates on a segment of the form $ax - (b-a)y = c$ with $c$ integer and $x \geq 0, y \leq 0$ as far down as possible. Note that the result of the action of generalized column regularization of a partition may not be a partition. Here is an example of the action of ${\rm colreg}_{3,5}$ on the partition $(3,2,1)$:

\

\ \  \  \   \   \   \   \    \ \ \ \ \ \ \ \ \ \ \  \
\begin{ytableau}
\ & & *(green)\\
 &  \\
 \\
\end{ytableau}  \   \  \    \    \ \    \  \   \ $\xrightarrow{\rm colreg_{3,5}}$ \  \  \   \  \  \  \  \  \ 
\begin{ytableau}
\ &  \\
\ & \\
\\
*(green)\\
\end{ytableau}

\

A box outside of a partition is called addable if the set union of this box and the partition is still a partition. A box in a partition is called removable if the set difference of the partition and this box is still a partition. For a given residue $i$ we define a sequence of addable and removable boxes with residue $i$ written from the Southwest to the Northeast, writing A or R for each box if it is addable or removable respectively. Then we cross out every word of the form RA from the sequence inductively. The first R from the left after all crossing-out is completed corresponds to a good box of residue $i$ in the partition with which we started (the notion of a a good box is due to Kleshchev, cf. \cite{K}).

To state Bezrukavnikov's conjecture, we need to define a new operation $\widetilde{M}^{\prime}_{a/b}$ on partitions of $n$ for a given term $a/b$ of the $n$-th Farey sequence. First we define the operation of concatenation $\cup$ of two partitions as follows: the sequence of rows of the new partition is the multiset union of the two sequences of rows of the old partitions. Next we write our partition as a concatenation $\mu = \nu \cup d \rho$ where each row of $\nu$ is not divisible by $d$ and each row of $d \rho$ is divisible by $d$, and is in fact equal to $d$ times the corresponding row of the partition $\rho$. Note that there is a unique way to write a partition in this way as a concatenation of two partitions, one of which has all rows not divisible by $d$ and the other one has all rows divisible by $d$. Finally we define $M^{\prime}_d(\mu) = \nu \cup d \rho^t$ and $\widetilde{M}^{\prime}_{a/b} = M_{b^{\prime}_k}^{\prime} \circ \dots \circ M_{b^{\prime}_1}^{\prime}$, where  $0<\frac{a^{\prime}_1}{b^{\prime}_1}<...<\frac{a^{\prime}_k}{b^{\prime}_k}<a/b$ be all the rational numbers with denominator at most $n$ between $0$ and $a/b$. 

A partition is called $d$-regular if none of its parts is repeated at least $d$ times. In the statement of Bezrukavnikov's conjecture the following extension of the Mullineux involution $M_d$ (where we also denote the extension by the same letter $M_d$) from the set of $d$-regular partitions to the set of all partitions is used. For a partition $\rho$ denote by $d * \rho$ the partition in which each part of $\rho$ is repeated $d$ times. Writing $\mu = \nu \cup d * \rho$ with $\nu$ regular, we define $M_d(\mu) = M_d(\nu) \cup d * \rho^t$. Now we are ready to state (see the definition of $\widetilde M_{a/b}$ in the introduction):

\begin{conjecture} {\rm (R. Bezrukavnikov)} For every positive integer $n$, every partition $\lambda$ of $n$, and every term of the Farey sequence $a/b$, the total number of boxes in all rows divisible by $b$ in the two partitions $\widetilde{M}_{a/b}(\lambda)$  and $\widetilde{M}^{\prime}_{a/b}(\lambda^t)$ is the same. 
\end{conjecture}

\begin{remark} Examples confirming this conjecture for $n \leq 5$ can be found in the appendix to \rm{\cite{DY}}.
\end{remark}

\section{Trivial representation}

Recall that we place Young diagrams in the fourth quadrant so that their horizontal and vertical boundaries are aligned with coordinate axes.

\begin{lemma} \label{lemma} For two neighboring fractions $\frac{a}{b}$ and $\frac{a^{\prime}}{b^{\prime}}$ in the Farey sequence for $n$ and two boxes $(x_i,y_i)$, $i=1,2$, in a Young diagram with $n$ boxes we cannot simultaneously have $a^{\prime}x_1 - (b^{\prime}-a^{\prime}) y_1 < a^{\prime} x_2-(b^{\prime}-a^{\prime}) y_2$ and $a x_1 - (b-a) y_1 > a x_2 - (b-a) y_2$.  
\end{lemma}

\begin{proof}

The two inequalities $ax_1-(b-a)y_1 > a x_2-(b-a) y_2$ and $a^{\prime}x_2-(b^{\prime}-a^{\prime})y_2 > a^{\prime}x_1-(b^{\prime}-a^{\prime})y_1$ which hold simultaneously lead to $(b-a)(y_2-y_1)>a(x_2-x_1)$ and $(b^{\prime}-a^{\prime})(y_1-y_2)>a^{\prime}(x_1-x_2)$. This in turn implies that $\frac{x_1-x_2}{y_1-y_2}$ is between $\frac{b-a}{a}$ and $\frac{b^{\prime}-a^{\prime}}{a^{\prime}}$. This means that $\frac{x_1-x_2+y_1-y_2}{y_1-y_2}$ is between $\frac{b}{a}$ and $\frac{b^{\prime}}{a^{\prime}}$. From this we conclude that $\frac{|y_1-y_2|}{|x_1+y_1-x_2-y_2|}$ is between $\frac{a}{b}$ and $\frac{a^{\prime}}{b^{\prime}}$. We have $|x_1+y_1-x_2-y_2| \leq |x_1-x_2|+|y_1-y_2| \leq n$ because $(x_1,y_1)$ and $(x_2,y_2)$ lie within a Young diagram with $n$ boxes total. This in turn contradicts the fact that $\frac{a}{b}$ and $\frac{a^{\prime}}{b^{\prime}}$ are neighboring terms of the Farey sequence and hence do not have any rational numbers with denominator at most $n$ between them.

\end{proof}

The proof of Theorem \ref{main} follows directly from the combination of the following two lemmas.

\begin{lemma} Let $\frac{a}{b}<\frac{a^{\prime}}{b^{\prime}}$ be the two neighboring terms in the $n$-th Farey sequence. Then $D^n_{a^{\prime},b^{\prime}}={\rm colreg}_{a,b}(D^n_{a,b})$.
\end{lemma} 

\begin{proof}

Suppose that the conclusion of the lemma is false. This means that there exist boxes $(x,y)$ in ${\rm colreg}_{a,b}(D^n_{a,b})$ and $(\widetilde x, \widetilde y)$ not in ${\rm colreg}_{a,b}(D^n_{a,b})$ such that $(\widetilde x, \widetilde y) <_{a^{\prime},b^{\prime}} (x,y)$.

Case 1. Assume that $a^{\prime} \widetilde x - (b^{\prime} - a^{\prime}) \widetilde y < a^{\prime} x - (b^{\prime} - a^{\prime}) y$. By Lemma \ref{lemma} we have $a \widetilde x - (b - a) \widetilde y \leq a x - (b - a) y$. If $a \widetilde x - (b - a) \widetilde y < a x - (b - a) y$ then since $(x,y)$ in ${\rm colreg}_{a,b}(D^n_{a,b})$ we have that $(\widetilde x, \widetilde y) \in D_{a,b}$, which is a contradiction. So $a \widetilde x - (b - a) \widetilde y = a x - (b - a) y$. Therefore $(x,y)$ and $(\widetilde x, \widetilde y)$ lie on the same line of slope $\frac{b-a}{a}$. Since $(x,y)$ in ${\rm colreg}_{a,b}(D^n_{a,b})$ and $(\widetilde x, \widetilde y)$ not in ${\rm colreg}_{a,b}(D^n_{a,b})$, we have that $(\widetilde x, \widetilde y)$ lies above $(x,y)$ i.e. $\widetilde y > y$. 

Therefore we obtain $a^{\prime} \widetilde x - (b^{\prime} - a^{\prime}) \widetilde y < a^{\prime} x - (b^{\prime} - a^{\prime}) y$, $a \widetilde x - (b - a) \widetilde y = a x - (b - a) y$, and $\widetilde y > y$. From the equality we get $\widetilde y - y = \frac{a}{b-a}(\widetilde x - x)$, hence $\widetilde x > x$. Substituting this into the inequality, we get $a^{\prime}(\widetilde x - x)<(b^{\prime} - a^{\prime})\frac{a}{b-a}(\widetilde x - x)$. This contradicts $\frac{a}{b}<\frac{a^{\prime}}{b^{\prime}}$. Therefore we obtained a contradiction in Case 1.  

Case 2. Assume that $a^{\prime} \widetilde x - (b^{\prime} - a^{\prime}) \widetilde y = a^{\prime} x - (b^{\prime} - a^{\prime}) y$ and $\widetilde y > y$. Since $(\widetilde x, \widetilde y)$ is outside of ${\rm colreg}_{a,b}(D^n_{a,b})$, we have that $a \widetilde x - (b - a) \widetilde y = a x - (b - a) y$. Therefore we have that $a^{\prime}(\widetilde x - x) = (b^{\prime} - a^{\prime})(\widetilde y - y)$ and $a(\widetilde x - x) \geq (b - a)(\widetilde y - y)$. From these two equalities we obtain $\frac{b^{\prime} - a^{\prime}}{a^{\prime}}(\widetilde y - y) \geq \frac{b - a}{a}(\widetilde y - y)$ which contradicts $\frac{a}{b}<\frac{a^{\prime}}{b^{\prime}}$. Therefore we obtained a contradiction in Case 2.   

\end{proof}

\begin{lemma} b-Mullineux transposed and the column regularization ${\rm colreg}_{a,b}$ have the same effect on $D^n_{a,b}$: we have $M_b^t(D^n_{a,b})={\rm colreg}_{a,b}(D^n_{a,b})$. 
\end{lemma}

\begin{proof} 

We define an order on pairs $(n, \frac{a}{b})$ where $\frac{a}{b}$ is a term of the Farey sequence for $n$ by $(n_1, \frac{a_1}{b_1})>(n_2, \frac{a_2}{b_2})$ if $n_1>n_2$ or $n_1=n_2$ and $\frac{a_1}{b_1}>\frac{a_2}{b_2}$. We prove this lemma by increasing induction on the pair $(n, \frac{a}{b})$.

Consider the box $P=D_{a,b}^n \setminus D_{a,b}^{n-1}$. We claim that this box is $b$-good by definition of a $b$-good box (see Section 2). We check the definition of a $b$-good box for $P$ as follows: we claim that all addable and removable boxes with the same residue as $P$ lie on the line through $P$ with slope $\frac{b-a}{a}$ above $P$ and from left to right all addable boxes go before all removable boxes, with the first removable box from the left in the sequence being the box $P$.

Let $P$ lie on the line $ax-(b-a)y=c$ in the diagram. 
 
Let us show that there are no addable or removable boxes in $D_{a,b}$ with the same $b$-residue as $P$ such that they lie above the line $ax-(b-a)y=c$. Suppose there is an addable or a removable box $(x_1,y_1)$ in $D_{a,b}$ with the same $b$-residue as $P$ above the line $ax-(b-a)y=c$. Then we have $ax_1-(b-a)y_1=c-kb$ for some positive integer $k$. Then $a(x_1+1)-(b-a)y_1=c-kb+a<c$ therefore the box $(x_1+1,y_1)$ lies in $D^n_{a,b}$ and so $(x_1,y_1)$ cannot be an addable or a removable box. Contradiction. 

Now let us show that every box on the line $ax-(b-a)y=c$ is either addable if it is to the left of $P$ or removable if it is to the right of $P$ or if it is $P$ itself. Let $Q$ be a box with coordinates $(x_2,y_2)$ on the line $ax-(b-a)y=c$. Then $a(x_2+1)-(b-a)y_2>c$ and $ax_2-(b-a)(y_2-1)>c$ hence the boxes immediately to the right of $Q$ and immediately below $Q$ are not in the diagram $D^n_{a,b}$. Hence $Q$ is addable if $Q$ is not in $D^n_{a,b}$ and removable if it is in $D^n_{a,b}$. Since all boxes on the line $ax-(b-a)y=c$ to the left of $P$ are not in $D^n_{a,b}$ and all boxes on the line $ax-(b-a)y=c$ to the right of $P$ and $P$ itself are in $D^n_{a,b}$, we obtain our conclusion.   

Similarly we find that the box $Y={\rm colreg}_{a,b}(D^n_{a,b}) \setminus {\rm colreg}_{a,b}(D^{n-1}_{a,b})$ is good in the partition which is the transpose of the partition ${\rm colreg}_{a,b}(D^n_{a,b})$ and has the opposite residue to the residue of $X$ in the partition $D^n_{a,b}$.

In Corollary 4.12 of \cite{BeO} (cf. p. 268 in \cite{FK}) it is shown that the result of the action of Mullineux on the difference of a partition and its good box with residue $i$ is the difference of Mullineux of the partition and its good box with residue $-i$. 

Since we showed above that $Y$ is the only addable box of residue $-i$ in the partiton which is the transpose of ${\rm colreg}_{a,b}(D^{n-1}_{a,b})$, by the above property of the action of Mullineux on a partition with a good box, $M_b(D^n_{a,b}) = M_b(D^n_{a,b} \cup X)=({\rm colreg}_{a,b}(D^{n-1}_{a,b}) \cup Y)^t = ({\rm colreg}_{a,b}(D^n_{a,b}))^t$.

\end{proof}

\begin{remark} As we increase n, to the left of the wall $a/b$ new boxes add to the Young diagram in the increasing order $>_{a,b}$.
\end{remark}

\begin{corollary} We deduce Theorem {\rm 3.5(2)} in {\rm \cite{DY}} which states that Mullineux transposed gives us generalized column regularization at each step of wall-crossing for the trivial representation.
\end{corollary}

\begin{corollary} Bezrukavnikov's conjecture for the trivial representation follows because it is easy to see that under both algorithms the partition stays regular all the time (under the algorithm of Section 2 the column partition never changes).
\end{corollary}

\section{Sign representation}

Let $n$ be a prime number and let $sgn_n$ be the partition of $n$ corresponding to the sign representation of the symmetric group $S_n$, i.e. its Young diagram is a column of height $n$. Let $r$ be a term in the $n$-th Farey sequence. Given positive integers $f_1<\dots<f_m$ and $g_1,\dots,g_m$, we will write $f_1^{g_1} f_2^{g_2} \dots f_m^{g_m}$ for the partition of the integer $f_1 g_1 + \dots + f_m g_m$ with $g_i$ parts equal to $f_i$. 

\begin{conjecture} \label{Mullsign} There is an increasing sequence of fractions $\frac{a_i}{b_i}$ in the Farey series such that for $\frac{a_i}{b_i}<r<\frac{a_{i+1}}{b_{i+1}}$ the partition $\widetilde{M}_r(sgn_n)$ of $n$ is of the form $z^t x^y$ where $x,y,z,t$ satisfy

(a) $z+y=b_i$

(b) $y+t+z-x=b_{i+1}$

\end{conjecture} 

\begin{remark} Together with the equation $xy+zt=n$, the above two equations form a system of three equations with four unknowns to be solved over the integers, which is equivalent to one linear Diophantine equation $xb_i-y(b_{i+1}-b_i)+b_i(b_{i+1}-b_i)=n$ with two unknowns.
\end{remark}

\begin{remark} One can begin proving the conjecture by noticing that Mullineux's original definition of the $b_i$-Mullineux involution uses the notion of $b_i$-rim of a partition. From Conjecture \ref{Mullsign} it follows that the $b_{i+1}$-rim of the partition $\widetilde{M}_r(sgn_n)$ is obtained from the usual rim of $\widetilde{M}_r(sgn_n)$ by removing the top row, and the $b_i$-rim of the transposed partition $\widetilde{M}_r(sgn_n)^t$ is obtained from the usual rim of $\widetilde{M}_r(sgn_n)^t$ by removing the top row. Furthermore, after removing the first $b_{i+1}$-rim, the Young diagram still has at most two possible lengths of rows.
\end{remark}

\section{Added good box conjecture}

The following was checked using a computer program for Young diagrams of size at most $40$:

\begin{conjecture}
If one Young diagram is obtained from another by adding one box to the first row in the interval $(0,1/n)$, these diagrams will differ by adding one good box to some row at any step of the combinatorial wall-crossing.
\end{conjecture}

\begin{remark} There are interpretations of this conjecture in terms of a crystal for the elliptic Hall algebra and in terms of components of the fixed point set of an element of finite order in the Hilbert scheme of points in the plane.
\end{remark}

\end{document}